\documentclass[12pt,a4paper]{article}
\pdfoutput=1
\usepackage{color}
\usepackage[usenames,dvipsnames,svgnames,table]{xcolor}
\definecolor{Light}{gray}{.85}
\usepackage{extarrows}
\usepackage{soul}
\usepackage[round]{natbib}
\bibpunct{(}{)}{;}{a}{,}{,}
\usepackage[hyperfootnotes=false,colorlinks=true]{hyperref}
\hypersetup{pdfcreator={PDFLaTeX mit dem hyperref Paket},%
  pdfsubject={Residues, higher arithmetic},%
  pdfauthor={Christian Siebeneicher},%
  pdftitle={Residues  :  The gateway to higher arithmetic I},%
  pdfkeywords={Carl  Friedrich   Gauss,  Disquisitiones,  Peter
  Gustav Lejeune Dirichlet,     Lectures on  number theory, 
  The  Shaping
    of    Arithmetic after C.~F.~Gauss's Disquisitiones
    Arithmeticae},%
 urlcolor=black,%
 colorlinks=true,%
  citecolor=black,%
  }
\makeatletter
\renewcommand\section{\@startsection{section}{1}{0mm}%
 {10pt}%
 {6pt}%
 {\bfseries}}
\renewcommand\subsection{\@startsection{subsection}{2}{0mm}%
 {10pt}%
 {6pt}%
 {\bfseries}} \makeatother

\date{ }
\begin{document}
\noindent{}%
\begin{center}
  \textbf{\large  Residues  :  The gateway to higher arithmetic I}\\[2ex]
  Recalling section one of Gauss's Disquisitiones
  Arithmeticae\\[2ex]
  \textsc{Christian
    Siebeneicher}
\end{center}

\begin{quote}
  \footnotesize  \textsc{Abstract}: Residues to  a given  modulus have
  been  introduced to  mathematics by  Carl Friedrich  Gauss  with the
  definition  of  congruence  in  the  `Disquisitiones  Arithmeticae'.
  Their  extraordinary properties provide  the basis  for a  change of
  paradigm in arithmetic.  By  restricting residues to remainders left
  over by divison Peter Gustav Lejeune Dirichlet --- Gauss's successor
  in G\"ottingen --- eliminated in his `Lectures on number theory' the
  fertile    concept   of   residues    and   attributed    with   the
  number--theoretic approach to residues for  more than one and a half
  centuries  to obscure  Gauss's  paradigm shift  in mathematics  from
  elementary to higher arithmetic.
\end{quote}

\section{The outset of Dirichlet's number theory in a nutshell}

\noindent{}%
The \emph{Lectures  on number theory} ---  posthumously published 1863
by  Richard Dedekind  --- open  as follows:  \emph{In this  section we
  treate a few arithmetic theorems,  which indeed may be found in most
  text books,  but which  are of such  fundamental importance  for our
  science        that         a        rigorous        proof        is
  necessary}.\footnote{\citet[1]{dirichlet63}.}  Section~1 starts with
the proof of following general theorem extending the range of familiar
elementary arithmetic:\footnote{\citet[2]{dirichlet63}.}
\begin{quote}\small When  one forms  the product of any two  out of
  three  numbers,  and then  multiplies  this  product  by the  third
  number, the resulting pro\-duct always has the same
  value, regardless of which two numbers are chosen first.\\
  Since  the product  is independent  of the  order of  the successive
  multiplications, it is called the product of the three numbers, and
  the latter  are called its factors, without regard  to order.
\end{quote}
Subsequently it is  shown that a similar theorem  holds for any system
$S$ of arbritrary positive integers $a,b,c, \dots{}$:
\begin{quote}\small
  The simplest way  a product can be formed from  these numbers is the
  following.  One takes  two numbers from S at  random and forms their
  product;  the System  $S'$ of  the latter  number and  the remaining
  numbers then  has one less number  than $S$.  By  again choosing two
  numbers at random from $S'$ and forming their product, one obtains a
  system $S''$  with the  two numbers fewer  than $S$.   Continuing in
  this way, one finally arrives at a single number, and the theorem to
  be  proved is that  \emph{the number  remaining at  the end  of this
    process  is   always  the  same,  no  matter   how  the  indiviual
    multiplications are ordered.}\\[.3ex]
To prove this we apply complete induction.
\end{quote}

The last paragraph of the introductory chapter to arithmetic captioned
with  ---  \emph{Looking  back}   ---  closes  with  an  authoritative
statement:\footnote{\citet[20]{dirichlet63}.}
\begin{quote}\small
  Here  we conclude  the series  of  theorems on  the divisibility  of
  numbers. But it  is worthwhiIe at this stage to  look back along the
  path these investigations have taken. It is now clear that the whole
  structure rests  on a single  foundation, namely the  algorithm for
  finding  the  greatest  common  divisor  of two  numbers.   All  the
  subsequent theorems, even when they  depend on the later concepts of
  relative  and   absolute  prime  numbers,  are   still  only  simple
  consequences of the this initial investigation, so one is entitled to
  make the following claim: any analogous theory, for which there is a
  similar algorithm  for the greatest  common divisor, must  also have
  consequences  analogous  to those  in  our  theory.   In fact,  such
  theories exist.
\end{quote}
\S.~17 of chapter  2 ---\emph{On the congruence of  numbers} --- opens
with     a    first     consequence    based     on     the    initial
investigation:\footnote{\citet[21]{dirichlet63}.}
  \begin{quote}\small
\emph{If $k$  is any positive integer,  then each integer  $a$, may be
  written in exactly one way in the form
\begin{displaymath}
  a= sk+r,
\end{displaymath}
where $s$ is an integer an $r$  is one of  the $k$ numbers}
\begin{displaymath}
0,  1,  2,  \dots{} ,  (k -1).  
\end{displaymath}
\hspace*{2em}\vdots

In what follows we shall say that $r$ is  the \emph{remainder} of the number
$a$ relative to the \emph{modulus} $k$.
\end{quote}

\S.~18 then  defines and treats to some  extent \emph{complete residue
  systems}  comprising  theorem 6  \emph{of  particular importance  in
  later investigations}:\footnote{\citet[24]{dirichlet63}.}
\begin{quote}\small
  If $a$ is relatively prime to  the modulus $k$ and if one replaces $x$
  in the linear expression $ax + b$  by the series of all $k$ terms in
  a complete system of  incongruent numbers, then the resulting values
  also form a complete system of incongruent numbers.
\end{quote}

The \emph{remainder} of a number relative to the \emph{modulus} $k$ is
used in \S.~19 for the  \emph{proof of the generalised Fermat theorem}
in   the  form  $a^{\varphi(k)}\equiv   1  \pmod   k$  which   may  be
\emph{expressed                in               words               as
  follows}:\footnote{\citet[25--26]{dirichlet63}.}
  \begin{quote}\small
If $a$ is relatativly prime to the positive integer $k$, and if  one
raises $a$ to the power  $\varphi(k)$, equal to the number of numbers among
\begin{displaymath}
   1,   2,   3, \dots{},       k
\end{displaymath}
relatively prime to $k$, then the result always leaves remainder 1 on division by $k$.
\end{quote}
Afterwards,  in \S.~20  \emph{another proof  of the  same  theorem} is
added with  the accentuation:  \emph{It is not  at all  superfluous to
  give another proof of this  important theorem which can be developed
  from the binomial theorem.}\footnote{\citet[27]{dirichlet63}.}

The explicit reference to \emph{the binomial theorem} calls to mind an
annotation  to  Fermat's  theorem   disclosing  in  section~3  of  the
\emph{Disquisitiones Arithmeticae}  a surprising detail  to the theory
of numbers:\footnote{\citet[32]{gauss86}.}
\begin{quote}\small
  This theorem is worthy of attention both because of its elegance and
  its great  usefulness.  It is usualIy called  Fermat's theorem after
  its  discoverer.  (See Fermat,  \emph{Opera Mathem.}   p.~131) Euler
  first published  a proof in his  dissertation entitled ``Theorematum
  quorundam  ad  numeros  primos  spectantium  demonstratio,''  (Comm.
  acad.           Petrop.,           8          [1736],          1741,
  141).\long\def\symbolfootnote[#1]#2{\begingroup%
    \def\thefootnote{\fnsymbol{footnote}}\footnote[#1]{#2}\endgroup}%
  \symbolfootnote[1]{Gauss's footnote: In a previous commentary (Comm.
    Acad.   Petrop.,  6   1723--33],  1738,  106  [``Observationes  de
    theoremate   quodam   Fermatiano   aliisque  ad   numeros   primos
    spectantibus'']), this great man  had not yet reached this result.
    In the  famous controversy between  Maupertius and K\"onig  on the
    principle of  the least action,  a controvery that led  to strange
    digressions,  K\"onig claimed  to  have at  hand  a manuscript  of
    Leibniz  containing  a demonstration  of  this  theorem very  Iike
    EuIer's (K\"onig, \emph{Appel au public},  p. l06). We do not wish
    to deny this testimony,  but certainly Leibniz never published his
    discovery.  Cf.  \emph{Hist.  de  l'Ac.de Prusse, 6} [1750], 1752,
    530. [The reference is ``Lettre de M. Euler \'a M. Merian (traduit
    du Latin)'' which was sent from Berlin Sept. 3, 1752.] } The proof
  is  based on the  expansion of  $(a +  1)^p$. From  the form  of the
  coefficients  it is easily  seen that  $(a +  1)^p -  a^p -  1$ will
  always be divisible by $p$ and consequently so will $(a + 1)^p- (a +
  1)$ whenever $a^p - a$ is divisible  by $p$.  Now since $1^p - 1$ is
  always divisible by  $p$ so also will  $2^p - 2$, $3^p -  3$, and in
  general $a^p - a$.  And since $p$ does not divide $a$, $a^p- 1$ will
  be divisible by $p$.  ---  \emph{That will suffice to to clarify the
    character of that method}.\footnote{The reference to the character
    of  the  two  different  methods  of  proof  of  Fermat's  theorem
    indicating  for the  first one  the involvement  of  the induction
    principle has been omitted in the English translation.  Because of
    its  tremendous  importance   for  the  underlying  principles  of
    arithmetic an English translation was inserted on the basis of the
    Latin   original.}   The  illustrious   Lambert  gave   a  similar
  demonstration in \emph{Nova Acta  erudit.}, 1769, p.~109.  But since
  the expansion of  a binomial power seemed quite  alien to the theory
  of  numbers,  Euler  gave  another  demonstration  that  appears  in
  \emph{Novi comm.  acad.  Petrop.}  T.   VIII p.~70, which is more in
  harmony with  what we have done  in the preceding  article.  We will
  offer  still  others later.   Here  we  will  add another  deduction
  similar to those of Euler~\dots{}
\end{quote}

The pointer to the preceding article which is said to be \emph{more in
  harmony  with what  we  have done}  leads  to take  note of  article
49:\footnote{\citet[30]{gauss86}.}
\begin{quote}\small
  49. \textsc{Theorem}.  \emph{If $p$ is  a prime that does not divide
    $a$, and  $a^t$ is the  lowest power of  $a$ that is  congruent to
    unity relative to the modulus $p$, the exponent $t$ will be either
    $= p- 1$ or be a factor of this number.}
\end{quote}

The  preliminary  note  preceding  the demonstration  of  the  theorem
states: \emph{We have  alread seen that $t$ either  $=p-1$ or $< p-1$.
  It remains  to show  that in the  latter case  $t$ will always  be a
  factor of $p  -1$}. That remarkable property stimulates to go further  back in section 3
thereby bringing to light details to the method of proof \emph{more in
  harmony } with arithmetic:\footnote{\citet[29]{gauss86}.}
\begin{quote}\small
  45. \textsc{Theorem}.  \emph{In any geometric progression $1,  a, aa, a^3$
  etc., outside of the first term~1, there is still another term $a^t$
  which is congruent to unity relative  to the modulus $p$ when $p$ is
  prime   relative  to   $a$,  and   the  exponent~$t$  is   $<  p$.}
\end{quote}
That theorem shows that the  residues are generated from the leading 1
by  repeated   multiplication  with  the  scale  factor   $a$  of  the
progression ---  an exceptional and far--reaching  phenomenon that has
not been accentuated Dirichlet's  lectures on number theory.  Moreover
section      three       with      caption      \emph{Residues      of
  Powers}\footnote{\citet[29--32]{gauss86}.}   starts  precisely  with
the general arithmetic truth providing an enlightening introduction to
the first  main theme of higher  arithmetic which has  been ignored in
Dirichlet's lectures on number theory:
\begin{center}
  \textsf{\scriptsize{}The  residues  of  the  terms  of  a  geometric
    progressions  which  begins   with  unity  constitute  a  periodic
    series}\footnote{Marginal notes of Waterhouse's revised version of
    Clarke's English  translation have been  transformed to catchlines
    as put  before in  Schering's edition of  the \emph{Disquisitiones
      Arithmeticae}  \citep{gauss63}. The  original edition  from 1801
    does  not  provide  that  guidance but  dispays  the  fundamental
    arithmetic truth on page~XIII of the table of contents.}
\end{center}

Before  coming to  the first  section  supplying the  basis of  higher
arithmetic  these remarkable  issues prompt  to consider  as  well the
preface of the \emph{Disquisitiones Arithmicae} with a splitting
of arithmetic:\footnote{\citet[xvii]{gauss86}.}
\begin{quote}\small
  However  what is  commonly called Arithmetic  hardly extends
  beyond  the art  of  enumerating and  calculating (i.e.   expressing
  numbers   by   suitabIe   symbols,   for  example   by   a   decimal
  representation,  and  carrying  out  arithmetic operations).   As  a
  resuIt it  seems proper to  call what had been  mentioned Elementary
  Arithmetic  and  to  distinguish  from it  Higher  Arithmelic  which
  properly  includes more general  inquiries concerning  integers.  We
  consider in the present volume only Higher Arithmetic.  --- Included
  under  the  heading ``Higher  Arithmctic''  are  those topics  which
  Euclid  treated  in  Book  VIl~ff.   with  the  elegance  and  rigor
  customary  among  the ancients;  but  that  restricts  to the  first
  commencements of this science.
\end{quote}

\section{Section one  of the \emph{Disquisitiones  Arithmicae}}

\begin{center}
  \textsc{\normalsize \textsc{Congruent Numbers in General}}\footnote{\citet[1]{gauss86}.}\\
  \textsf{\scriptsize{}Congruent   numbers,   modules,  residues   and
    nonresidues}
\end{center}

\begin{quote}\small
1.   If   a  number   $a$  divides  the difference of
  the numbers $b$ and $c$, $b$  and $c$ are said to be \emph{congruent
    relative to  a}; if not,  \emph{noncongruent}.  The number  $a$ is
  called  the  \emph{modulus}.   If   the  numbers  $b$  and  $c$  are
  congruent, each of them is called a \emph{residue} of the other.  If
  they are noncongruent they are called \emph{nonresidues}.

  The numbers involved must  be positive or
  negative  integers,\long\def\symbolfootnote[#1]#2{\begingroup%
  \def\thefootnote{\fnsymbol{footnote}}\footnote[#1]{#2}\endgroup}%
\symbolfootnote[1]{Gauss's footnote:  The modulus
    must  obviously be  taken  absolutely, i.e.   without sign.}   not
  fractions.  For  example, $- 9$ and  $+ 16$ are  congruent modulo 5;
  $-7$ is  a residue of $+15$  modulo 11, but  a nonresidue modulo~3. 

Since every  number divides  zero, it follows  that we can  regard any
number as congruent to itself relative to any modulus.
\end{quote}

\definecolor{gray90}{gray}{0.80}       \definecolor{gray85}{gray}{0.75}
\definecolor{gray80}{gray}{0.80}       \definecolor{gray75}{gray}{0.75}
\definecolor{gray70}{gray}{0.70}       \definecolor{gray65}{gray}{0.65}
\definecolor{gray60}{gray}{0.60}       \definecolor{gray55}{gray}{0.55}
\definecolor{gray50}{gray}{0.50}       \definecolor{gray45}{gray}{0.45}
\setlength{\fboxsep}{1.5pt}  The numerical  presetting  $-9 \equiv  16
\pmod 5$ suggestst  to mark with the numbers $-9$  and $16$ two points
on the integer line and to consider with these the entire intervall of
residues  modulo 5  located in  between.  Coloring  congruent residues
with the same color makes visible a geometric structure on the integer
line resulting for  the modulus 5 in consequence  of the definition in
article~1:
\begin{center}\scriptsize%
  \colorbox{gray65}{$-9$}%
  \colorbox{gray90}{$-8$}\colorbox{gray90}{$-7$}%
  \colorbox{gray65}{$-6$}\colorbox{gray45}{$-5$}\colorbox{gray65}{$-4$}%
  \colorbox{gray90}{$-3$}\colorbox{gray90}{$-2$}%
  \colorbox{gray65}{$-1$}%
  \colorbox{gray45}{$\phantom{\rule[-.2ex]{3pt}{1ex}}0\phantom{\rule[-.2ex]{3pt}{1ex}}$}%
\colorbox{gray65}{$\phantom{\rule[-.2ex]{3pt}{1ex}}1\phantom{\rule[-.2ex]{3pt}{1ex}}$}%
 \colorbox{gray90}{$\phantom{\rule[-.2ex]{3pt}{1ex}}2\phantom{\rule[-.2ex]{3pt}{1ex}}$}%
\colorbox{gray90}{$\phantom{\rule[-.2ex]{3pt}{1ex}}3\phantom{\rule[-.2ex]{3pt}{1ex}}$}%
\colorbox{gray65}{$\phantom{\rule[-.2ex]{3pt}{1ex}}4\phantom{\rule[-.2ex]{3pt}{1ex}}$}%
  \colorbox{gray45}{$\phantom{\rule[-.2ex]{3pt}{1ex}}5\phantom{\rule[-.2ex]{3pt}{1ex}}$}%
  \colorbox{gray65}{$\phantom{\rule[-.2ex]{3pt}{1ex}}6\phantom{\rule[-.2ex]{3pt}{1ex}}$}%
  \colorbox{gray90}{$\phantom{\rule[-.2ex]{3pt}{1ex}}7\phantom{\rule[-.2ex]{3pt}{1ex}}$}%
  \colorbox{gray90}{$\phantom{\rule[-.2ex]{3pt}{1ex}}8\phantom{\rule[-.2ex]{3pt}{1ex}}$}%
  \colorbox{gray65}{$\phantom{\rule[-.2ex]{3pt}{1ex}}9\phantom{\rule[-.2ex]{3pt}{1ex}}$}%
  \colorbox{gray45}{$\phantom{\rule[-.2ex]{1pt}{1ex}}10\phantom{\rule[-.2ex]{1pt}{1ex}}$}%
\colorbox{gray65}{$\phantom{\rule[-.2ex]{1pt}{1ex}}11\phantom{\rule[-.2ex]{1pt}{1ex}}$}%
\colorbox{gray90}{$\phantom{\rule[-.2ex]{1pt}{1ex}}12\phantom{\rule[-.2ex]{1pt}{1ex}}$}%
\colorbox{gray90}{$\phantom{\rule[-.2ex]{1pt}{1ex}}13\phantom{\rule[-.2ex]{1pt}{1ex}}$}%
\colorbox{gray65}{$\phantom{\rule[-.2ex]{1pt}{1ex}}14\phantom{\rule[-.2ex]{1pt}{1ex}}$}%
  \colorbox{gray45}{$\phantom{\rule[-.2ex]{1pt}{1ex}}15\phantom{\rule[-.2ex]{1pt}{1ex}}$}%
  \colorbox{gray65}{$\phantom{\rule[-.2ex]{1pt}{1ex}}16\phantom{\rule[-.2ex]{1pt}{1ex}}$}%
\end{center}

The sum of the  residues 9 and 7 may be determined  by starting from 9
and stepping forward  7 steps. Surprisingly however, the  result 16 is
congruent as  well to any of  $-9, -4, 1,  6$ and 11 and  already that
peculiarity  indicates that  higher arithmetic  obeys rules  which are
considerably different  from those governing  elementary arithmetic as
taught to school children.

\begin{quote}\small
  2.   Given $a$, all  the residues  modulo $m$  are contained  in the
  formula  $a +  km$ where  $k$  is any  integer.  The  easier of  the
  propositions that we state below  follow at once from this, but with
  equal ease  they can be proved  directly, \emph{as will  be clear to
    the reader}.\footnote{\citet[1]{gauss86}.   The emphasized text is
    missing in  Waterhouse's revised edition  and has been  taken from
    Clarke's translation from 1966.}

  \long\def\symbolfootnote[#1]#2{\begingroup%
    \def\thefootnote{\fnsymbol{footnote}}\footnote[#1]{#2}\endgroup}
  Henceforth  we shall  designate  congruence by  the symbol  $\equiv$
  joining to it in parentheses the  modulus when it is necessary to do
  so;      e.g.      $-7\equiv      15      \pmod{11},\;     -16\equiv
  9\pmod{5}$.\symbolfootnote[7]{Gauss's footnote: We have adopted this
    symbol  bccause of  the analogy  between equality  and congruence.
    For the same reason Legendre, in the treatise which we shall often
    have  occasion  to cite,  used  the  same  sign for  equality  and
    congruence. To avoid ambiguity we have made a distinction.}
\end{quote}

The  formula $a +  km$ from  article~2 then  implies that  the pattern
shown  for the modulus  5 is  periodic and  repeats after  five steps.
Moreover, applying  the symbol $\equiv$  also to the sum  $9+7=16$ and
its different  manifestations as $9+7=16\equiv  -9 \equiv -4  \equiv 1
\equiv  6 \equiv  11$ make  visible a  phenomenon typical  for higher
arithmetic which is inconceivable in it's elementary shaping.

Article 3 then provides a straightforward means to assign to any given
integer  $A$ its uniquely  determined residue  modulo $m$:\footnote{In
  order  to  avoid  a line  break  inside  the  crucial datum  of  $m$
  successive   residues  the   layout  of   the   English  translation
  \citep[1--2]{gauss86}  has been changed  to Schering's  reissue from
  1863 \citep[10]{gauss63} of the original latin text from 1801.  Like
  changes of  the English translation will occur  also later \emph{for
    the ease of order and summary} \citep[vi]{gauss86}.}

\begin{quote}\small
  3.  \textsc{Theorem}.   \emph{Let  m  successive
  integers
  \begin{displaymath}
    a,\, a  + 1,\, a + 2,\,  \dots{}\,, a + m - 1  
  \end{displaymath}
  and another integer
  $A$ be  given; then  one, and  only one, of  these integers  will be
  congruent to $A$ modulo $m$.}

If $\frac{a-A}{m}$ is an integer  then $a\equiv A$; if it  is a fraction,
let $k$  be the next  larger integer (or  if it is negative,  the next
\emph{smaller} integer not regarding sign). $A + km$ wilI fall between
$a$ and  $a +m$  and will  be the desired  number.  Evidently  all the
quotients
\begin{displaymath}\small
\tfrac{a -A}{m}\,, \quad \tfrac{a  + 1 - A}{m}\,, \quad \tfrac{a + 2
  - A}{m}\,, \quad \mathrm{etc}.   
\end{displaymath}
 lie between
$k - 1$ and $k + 1,$ so only one of them can be an integer.
\end{quote}

Clearly  that  theorem makes  obsolet  Dirichlet's  algorithm for  the
greatest  common divisor.   Moreover, an  application of  that theorem
provides  on the  basis of  \emph{least  residues} the  first step  to
establish on the integer line an arithmetic of residues.
\begin{center}
  \textsf{\scriptsize{}  Least  residues}
\end{center}
\begin{quote}\small
  4.   Each  number  therefore  will  have a  residue  in  the  series
  \mbox{$0,1,2,\dots{}   m-    1$}   and   in    the   series   $0,-1,
  -2,\dots{} -(m-1)$.  We will call  these the  \emph{least residues},
  and it is  obvious that unless 0 is a residue,  they always occur in
  pairs,  one \emph{positive}  and one  \emph{negative}.  If  they are
  unequal in  magnitude one will be  $< m/2$ ; otherwise  each will $=
  m/2$ disregarding sign. Thus each  number has a residue which is not
  larger   than   half  the   modulus.    It   will   be  called   the
  \emph{absolutely least} residue.\\
For example, relative to the modulus 5, $-13$ has 2 as least
positive residue. It is also the absolutely least residue, whereas
$-3$  is the least negative residue. Relative to the modulus 7, $+5$ is its
own least positive residue; $-2$ is the least negative residue and the
absolutely least.                                  
\end{quote}
The  definition  of  \emph{least  positive residues}  and  \emph{least
  negative  residues}  permits  to  modify  the residues  of  any  two
integers $a$ and $b$ and their sum  $a+b$ in a way such that all three
are  contained in  one  and the  same  intervall of  length $m$.   The
following example shows for the modulus  5 how the sum $9+7=16$ may be
determined in this manner: $9+7 \equiv 4 +2 \equiv -1 + 2 = 1$.

Since for any given modulus $m$ only the finitely many residues in the
set $\{0, 1,  \dots, m-1\}$ are concerned a  general argument directed
to all integers as in Dirichlet's improvements of elemenary arithmetic
provided  by the \emph{few  arithmetic theorems,  which indeed  may be
  found  in  most  text  books,  but which  are  of  such  fundamental
  importance for our  science that a rigorous proof  is necessary} are
not necessary for higher arithmetic.

Any given arithmetic expression of residues may be treated for element
by element and the independence of the value of a product of the order
of the individual  factors for example may be  proved by the verifying
finitely many  cases.  That  renders unnecessary to  involve schematic
\emph{complete   induction}  directed   to   \so{\emph{all}}  positive
integers  which  therefore may seem \emph{quite alien} to arithmetic.
\begin{center}
  \textsf{\scriptsize{} Elementary propositions regarding congruences}
\end{center}
\begin{quote}\small
5.  Having  established  these  concepts,  let  us  establish  the
properties  that foIlow  from them.

\emph{Numbers that  are congruent relative to a  composite modulus are
  also congruent relative to any divisor of the modulus.}

\emph{lf many numbers are congruent to the same
number relative to  the  same modulus, they are congruent to one
another (relative to the same modulus).}

This identity of moduli is to be understood also in what follows.

\emph{Congruent numbers  have the same least  residues, noncongruent numbers
have different least residues.}

6. \emph{Given  the numbers $A, B,  C,$ etc. and other  numbers $a, b,
  c,$  etc.    congruent  to  each  other  relative   to  any  modulus
  whatsoever, i.e.  $A \equiv a, B \equiv  b,$ etc., then $A + B + C +
  etc. \equiv a + b + c + etc.$}

\emph{lf $A\equiv a, B\equiv b$, then $A - B\equiv a - b$.}

7. \emph{If $A \equiv a$ then also $kA\equiv ka$.}

If $k$ is a positive number, then this is onIy a particuIar case of
the preceeding article (art. 6) letting $A=B=C$ etc., $a=b=c$ etc.
If $k$ is negative, $-k$ wilI be positive. Thus $-kA\equiv-ka$ and so
$kA\equiv ka$.

 lf $A\equiv a, B \equiv b$ then $AB \equiv ab$ because $AB \equiv Ab
 \equiv ba$.

 8.  \emph{Given any  numbers whatsoever  $A,  B, C,$  etc. and  other
   numbers $a, b, c,$ etc.  which are congruent to them, i.e. $A\equiv
   a, B \equiv b,$ etc., the  products of each will be congruent, i.e.
   $ABC$ etc. $\equiv abc$ etc.}

From the prcceding articIe $AB \equiv  ab$  and for the same reason
$ABC\equiv  abc$ and any number of factors may be adjoined.

If all  the numbers $A, B, C,$  etc. and the corresponding  $a, b, c,$
etc.  are assumed equal, then the foIlowing theorem holds: If $A\equiv
a$ and $k$ is a postive integer $A^k\equiv a^k$.

9. Let $X$ be an algebraic function with undetermined $x$ of the form
\begin{displaymath}
  Ax^a + Bx^b + Cx^c + \; \mathrm{etc.}
\end{displaymath}
where $A,B,C,$  etc. are any integers; $a,b,c,$  etc.  are nonnegative
integers.  Then if $x$ is given values which are congruent relative to
some modulus,  the resulting values of  the function $X$  will also be
congruent.

Let  $f,  g$ be  congruent  values  of  $x$. Then  from  the
preceding article  $f^a\equiv g^a$ and $Af^a\equiv Ag^a$ and in the
same way  $Bf^b\equiv Bg^b$ etc.Thus
\begin{displaymath}
  Af^a+Bf^b+Cf^c + \;\mathrm{etc.} \equiv  Ag^a+Bg^b+Cg^c +
  \;\mathrm{etc.}  \quad \textsc{q.e.d.}
\end{displaymath}
It is easy to understand how this theorem can be extended to
functions of many undetermined variables.
\end{quote}

Apparently  Gauss's  paradigm  shift  regarding mathematics  with  the
transition  from  elementary  to higher  arithmetic\footnote{Evicently
  Thomas  Kuhn has  not yet  been able  to point  to  that significant
  example  in  \emph{The  Structure  of Scientific  Revolutions}  from
  1962.}  has not yet been implemented in mathematics and the primordial
importance of  residues as autonomous entities  of arithmetic suitable
to     reveal    unsuspected     facets    of     the     God    given
integers\footnote{According  to  Heinrich \citet[19]{weber91}  Leopold
  Kronecker (1823--1891) once proclaimed: ``God made the integers; all
  else  is the  work  of man'';  see  also \emph{Namenverzeichnis}  in
  Helmut  \citet[1,   467]{hasse50}.}   still  looks   forward  to  be
recognized.

That  fact  manifests  itself  in  a review  from  1995  \emph{On  the
  arithmetic methods of mathematicians  of the seventeenth century} in
which            Igor            Rostislavovich            Shafarevich
wonders:\footnote{\citet[1]{luzin95}.}
\begin{quote}\small 
  It remains a mystery why such simple objects as the integers require
  for  their  understanding  practically  the  whole  machinery  which
  mathematicians are  able to create.  But this  mystery is completely
  analogous  to  that enigmatic  parallelism  between mathematics  and
  physics discussed by  a number of men of science.   In any case both
  phenomena  are too  universal to  elucidate them  using  an improper
  stage of development.  Apparently here we have to do with a cardinal
  phenomenon:  the human  way of  thinking  and the  structure of  the
  Cosmos are parallel to each other.
\end{quote}

The obvious question  if the mystery results all  alone from the human
way  of  thinking  with regard  to  the  integers  is not  taken  into
consideration.



\bibliographystyle{apalike}
\bibliography{Littlebook}

\newcommand{\noopsort}[1]{} \newcommand{\printfirst}[2]{#1}
  \newcommand{\singleletter}[1]{#1} \newcommand{\switchargs}[2]{#2#1}
\begin{thebibliography}{}

\bibitem[Gau{\ss}, 1863]{gauss63}
Gau{\ss}, C.~F. (1863).
\newblock {\em Disquisitiones Arithmeticae, \em Werke : Herausgegeben von der
  k{\"o}niglichen Gesellschaft der Wissenschaften zu G{\"o}ttingen, Editor: E.
  Schering}.
\newblock Universit{\"a}ts--Druckerei, G{\"o}ttingen.

\bibitem[Gauss, 1986]{gauss86}
Gauss, C.~F. (1986).
\newblock {\em Disquisitiones Arithmeticae, \em translated by Arthur A. Clarke,
  Yale University Press, New Haven, 1966, revised by W. Waterhouse, et. al.}
\newblock Springer, New York, Berlin, Heidelberg, Tokyo.

\bibitem[Hasse, 1950]{hasse50}
Hasse, H. (1950).
\newblock {\em Vorlesungen {\"u}ber Zahlentheorie}.
\newblock Springer, Berlin.

\bibitem[{Lejeune Dirichlet}, 1999]{dirichlet63}
{Lejeune Dirichlet}, J. P.~G. (1999).
\newblock {\em Vorlesungen {\"u}ber Zahlentheorie, \em ed. R.~Dedekind.
  Braunschweig: Vieweg. English transl. J. Stillwell. History of Mathematics
  Sources 16.}
\newblock American Mathematical Society -- London Mathematical Society,
  Providence -- London.

\bibitem[Luzin, 1995]{luzin95}
Luzin, N.~N. (1995).
\newblock On the arithmetic methods of mathematicians of the seventeenth
  century (preface of {L. A. T}er--{M}ika{\`e}lyan's book).
\newblock {\em Mathematical reviews \em -- MR1266620 (95g:01013)}.

\bibitem[Weber, 1892]{weber91}
Weber, H. (1891--1892).
\newblock Leopold {K}ronecker ({N}achruf).
\newblock {\em Jahresbericht der Deutschen Mathematiker--Vereinigung}, 2:5--31.

\end{thebibliography}
\end{document}